\documentclass[]{article}

\usepackage{amssymb} 
\usepackage{amsmath} 
\usepackage{latexsym}
\usepackage{graphicx}

\usepackage{tikz}
\usetikzlibrary{calc}
\usepackage{pgf,tikz}
\usetikzlibrary{arrows}
\usepackage{hyperref}
\usepackage{amscd}
\usepackage{relsize}

\usepackage{tikz-cd}

\usepackage{ color, amsmath, amssymb, amsfonts, amscd, graphicx,latexsym, hyperref}
\usepackage{caption, subcaption}
\usepackage[pdf]{pstricks}

\usepackage[affil-it]{authblk}

\usepackage{amsthm}

\newtheorem{theorem}{Theorem}
\newtheorem*{utheorem}{Theorem}
\newtheorem{proposition}{Proposition}

\newtheorem*{uremark}{Remark}
\newtheorem{lemma}{Lemma}

\newtheorem*{udefinition}{Definition}
\newtheorem{corollary}{Corollary}

\title{Flat Surfaces with Finite Holonomy Groups}
\author{ \.{I}smail Sa\u{g}lam \thanks{Electronic address: \texttt{isaglamtrfr@gmail.com}  }}
\affil{Adana Science and Technology University  }

\date{}
\begin{document}

\maketitle

\begin{abstract}
We prove that flow of a generic geodesic  on a flat  surface with finite holonomy group is ergodic. We use this result to prove that  flows of generic billiards  on certain flat surfaces with boundary are
also ergodic.
\end{abstract}
\tableofcontents
\section{Introduction}
$\:$ $\:$ A  surface is called flat if it has a flat metric having finitely many singularities of conical type. These surfaces may not be orientable and they may have boundary components. Ergodicity of a generic  geodesic (or billiard trajectory), existence of closed geodesics (or  billiard trajectories) are the main problems studied in dynamics of flat surfaces and
billiards.

If a flat surface is orientable, closed and its holonomy group is trivial, then it is called very flat.
Above problems  admit satisfactory and complete answers for very flat surfaces. See \cite{Kerckhoff- Masur} for ergodicity of geodesic flows and \cite{Eskin-Asymptotic}, \cite{Masur-Growth}, \cite{Masur-Lower bound} for counting problem of closed geodesics. Also, see \cite{zorich}, \cite{ergodic-translation}, \cite{
	rational-billiards} for introductory presentations about ergodic theory of very flat surfaces.
Note that ergodicity of the geodesic flow immediately implies existence of asymptotic cycles. See \cite{Schwartzmann-Asiymtotic} for a definition of asymptotic cycle.

If we are given a rational polygon, we can  obtain a very flat surface by gluing a finite number of its copies along their edges. See \cite{Zemlyakov-topological-transitivity},\cite{coding-problem}, \cite{Kerckhoff- Masur}. Since directional flow at almost all directions of a very flat surface is ergodic, this implies that flow of a generic billiard in such a polygon is ergodic. See \cite{rational-billiards}, \cite{gutkin} for much information about ergodicity of billiard flows in rational polygons.

Our aim is to prove that flow of a generic geodesic on a closed, orientable flat surface  with finite holonomy group is ergodic. We call these surfaces really flat. Note that we say flow of a geodesic $g_t$ on such a  surface $S$ is ergodic if at each region of $S$ it spends a time proportional to the volume of the region.  Equivalently, this means that for each continuous $f$ on $S$, the following limit exists:

\begin{align*}
\lim_{K \to \infty}\frac{1}{K}\int_{0}^{K} f(g_t) dt \ =
\frac{1}{A(S)}\int_S f d\mu,
\end{align*}

\noindent where  $A(S)$ is the area of $S$.

 We  show that these surfaces  can be covered by
a very flat  surfaces, where the covering is branched, Galois and respects flat  metrics on the surfaces. Indeed, this
immediately implies that flow of a generic geodesic on a really flat surface is ergodic. Finally, we use these results to prove ergodicity of flows of generic billiards on certain flat surfaces with boundary. We call these surfaces also really flat. Note that these surfaces may be non-orientable.

\paragraph{Convention:} Surfaces that we consider and maps between them are  of class $C^{\infty}$. They are also compact. They may have boundary components and be non-orientable.

\section{Flat coverings}
$\:$ $\:$  In this section, we introduce a notion of covering for flat surfaces. From the topological point of view, such a covering is nothing else than a branched covering. But, we also  require them  to respect flat metrics.

Given a branched covering  $\psi: S^* \rightarrow S$ between two surfaces, we denote set of its branched  points and ramification points by $\frak{b}$ and $\frak{r}$, respectively. 

Let  $S^*$ and $S$ be two flat surfaces.

\begin{udefinition}
	A map $\psi:\ S^* \rightarrow S $ is called a flat (covering) map, if it is a
	branched covering, and
	$$\psi\restriction_{S^*-\psi^{-1}(\frak{b})}: \  S^*-\psi^{-1}(\frak{b}) \rightarrow S-\frak{b}$$
	is a local isometry.
\end{udefinition}

\noindent We describe two ways to obtain flat  maps.

\begin{itemize}
	\item
	Let  $G$ be a finite isometry group of  $S^*$ so that each element except identity fixes a finite number of points.
	Consider $C^{\infty}$ (topological) orbifold $S^*/G$.  $\psi: S^* \rightarrow S^*/G $ induces a flat metric on $S^*/G$, and $\psi$ is a flat map with respect to the metric on
	$S^*$ and the induced metric on $S^*/G$.
	
	\item
	Let $\psi: S' \rightarrow S$ be a branched cover. Observe that $\psi$ induces
	a flat metric on $S'$ even if it is not Galois, and this makes $\psi$
	a flat  map.
\end{itemize}

\begin{uremark}
	
	If $\psi: S^* \rightarrow S$ is a flat map and $x \in S^*$ is a point with  ramification index  $m$, then we have
	$$\theta(x)=m\theta(\phi(x)),$$
	where $\theta(x)$ and $\theta(\phi(x))$ are the angles at $x$ and $\phi(x)$, respectively.
\end{uremark}
\section{Really flat  surfaces}
$\:$ $\:$ $\:$ \ In this section, we introduce a new family of flat surfaces.  We do not exclude non-oriantable surfaces and surfaces with boundary from our discussion.

Recall that an orientable flat surface without boundary    is called very flat if it  has trivial holonomy group. Surfaces defined below are natural 
generalizations of the very flat surfaces.

\begin{udefinition}

	A flat surface $S$ is called really flat if it has the following properties:
	
	\begin{enumerate}
		\item
		Its holonomy group is finite.
		
		\item
		For any two non-singular points $x,y$ on the boundary of $S$,
	 for any curve $L$ joining them and for any non-zero vector $v\in T_x(S)$ we have
		$$\frak{t}_y(v') - \frak{t}_x(v) \ \text{is} \ \pi-\text{rational,}$$

	\end{enumerate}
 where $v'$ is parallel transport of $v$ through $L$, $\frak{t}_x(v)$ is the angle between $v \in T_x(S)$ and its boundary component, $\frak{t}_y(v')$ is the angle between $v'$ and the boundary component that  $y$ is in. 
\end{udefinition}
	
Observe that last item of the above definition  makes sense: rationality of the change in angles
does not depend on orientation of the boundary component.

\begin{figure}
	\begin{center}
		\includegraphics[scale=0.45]{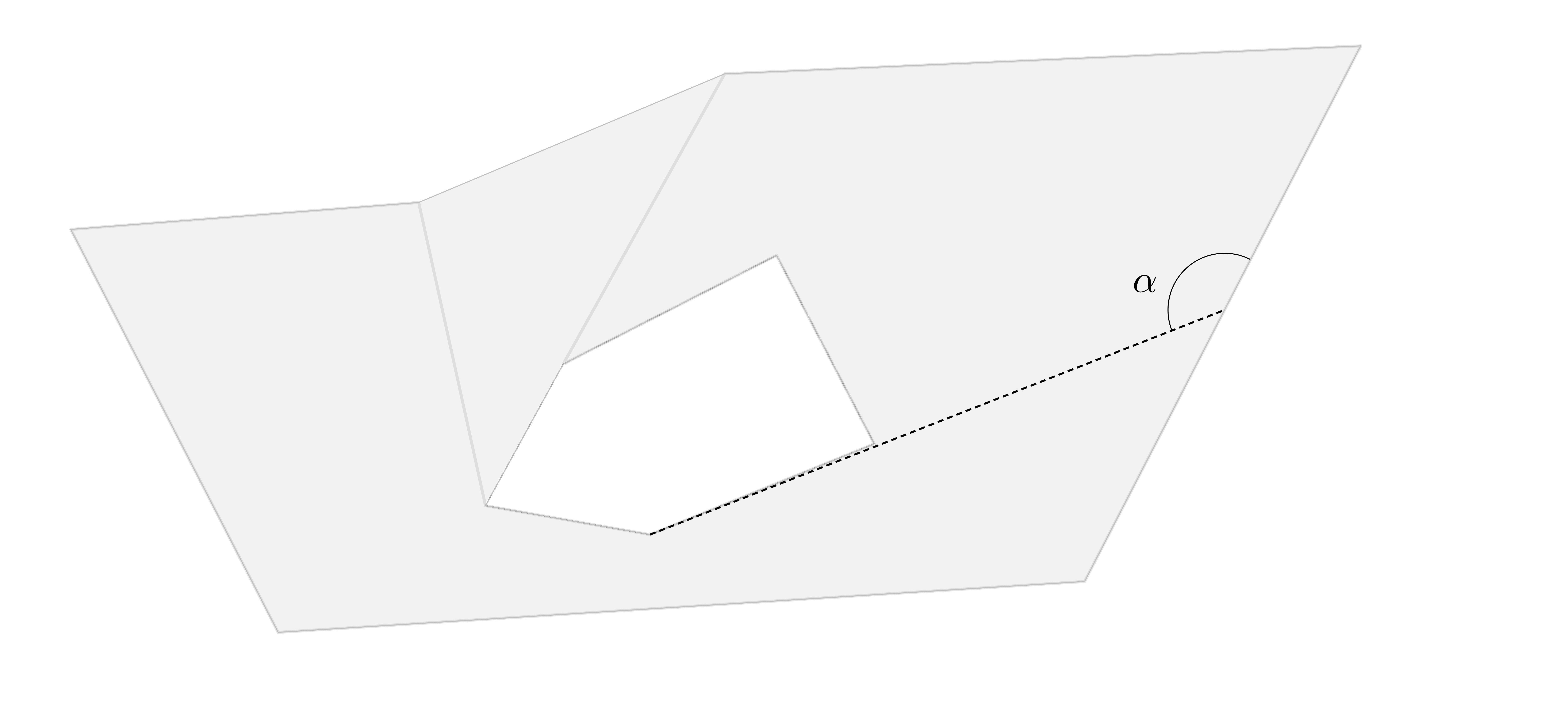}
	\end{center}
	\caption{A polygon in the plane. It is really flat if and only if $\alpha$ and angles at vertices of the polygon are $\pi$-rational.}
	\label{rationalpolygon}
\end{figure}

\begin{uremark}
\begin{enumerate}
	\item 
	If the surface does not have boundary component then it is really flat if and only if its holonomy group is finite.
	\item
	A planer polygon is rational if and only it is really flat. See Figure \ref{rationalpolygon}.
	\item
	Angles between boundary components of a really flat surface are $\pi$-rational.
	\end{enumerate}
\end{uremark}

\paragraph{Double of a flat surface} Let $S$ be a flat surface with at least one boundary component. Take another
copy of $S$, call it $S'$. For each  $x \in S$, let $x'$  be the corresponding point in $S'$. Glue $S$ and $S'$ so that each boundary point $x \in S$ comes together with $x'$. We will call resulting flat surface \textit{double} of $S$ and denote as $S^D$.

\begin{uremark}

	$S$ is really flat if and only if $S^D$ is really flat.
	
\end{uremark}

\section{Holonomy representation}

$\:$ $\:$  In this section, we assume that flat surfaces are   closed, oriented and 
flat maps are  orientation preserving. Our aim is to relate holonomy representations and flat maps.

Let $S$ be a flat surface, $\frak{s}$ be the set of  singular points of $S$ and $\frak{c}$ be a finite subset of $S$ that contains $\frak{s}$. Let $y \in S-\frak{c}$. We denote holonomy representation as
$$hol: \Pi_1(S-\frak{c},y) \rightarrow \mathbb{T}_y=\mathbb{T}=\{z \in \mathbb{C}: \lvert z \rvert =1\},$$
where $\mathbb{T}_y$ is the rotations of the  unit circle in the tangent space at $y$ of $S$, and  $\Pi_1(S-\frak{c},y)$ is the fundamental group
of $S-\frak{c}$ based at $y$.

\begin{udefinition}
	A curve on $S$ is called polygonal if it consists of finitely many geodesic segments and  does not pass through any point in $\frak{s}$.
\end{udefinition}
\begin{figure}
	\begin{center}
		\includegraphics[scale=0.45]{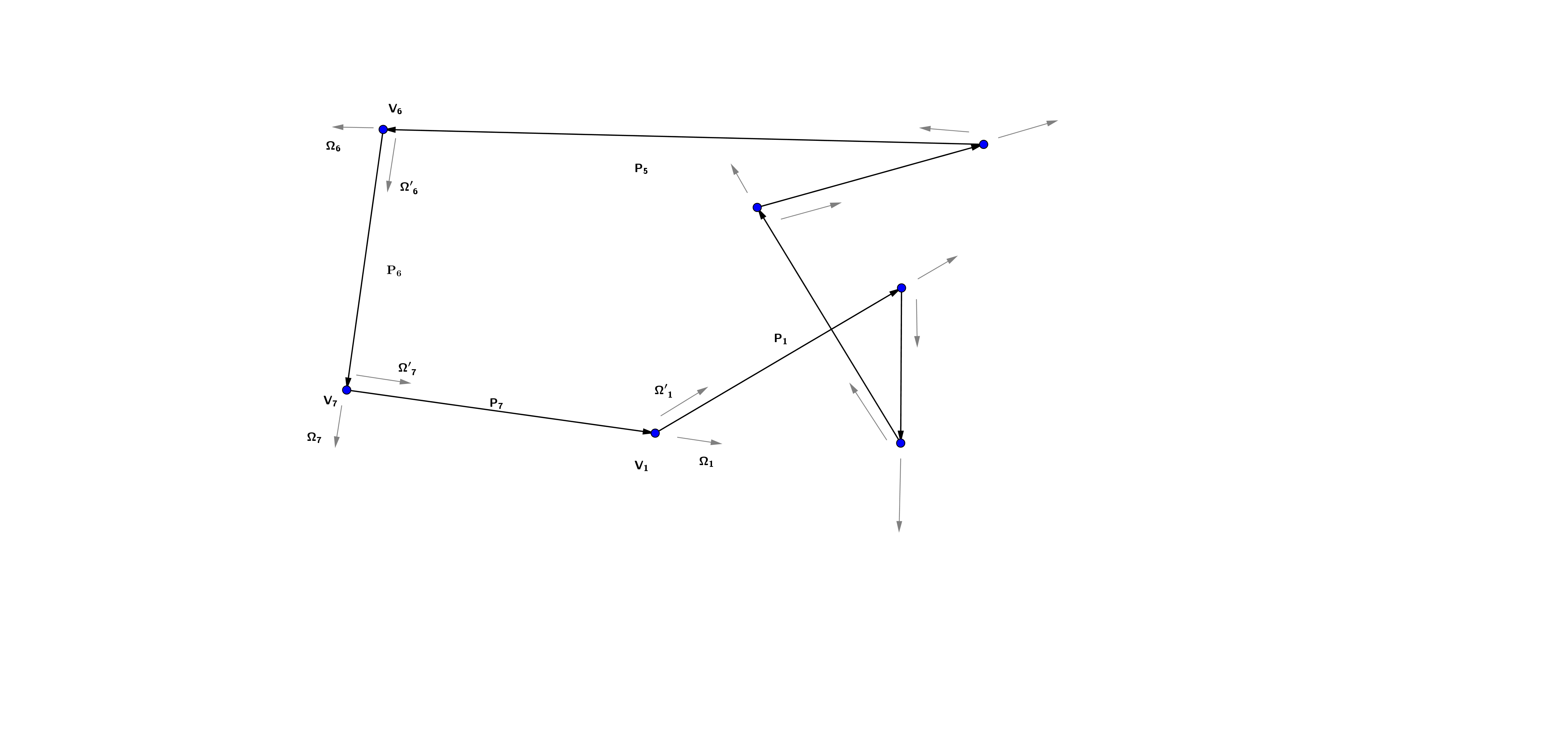}
	\end{center}
	\caption{A polygonal loop with $7$ vertices.}
	\label{polygonal-loop}
\end{figure}
Note that we do not assume  a polygonal curve is not self-intersecting.
Even it is possible that a geodesic segment of a polygonal curve is self-intersecting.

Observe that we can define (oriented) angle between two non-zero vectors $\beta_1, \beta_2$ at a given point of $S$. We will denote the angle
by $\theta(\beta_1,\beta_2)$.

Let $P$ be a polygonal  loop  based at $x \in S$. Let $n$ be the number of its vertices, and denote these vertices by $v_1, \dots, v_n$  so that ordering is compatible with orientation of $P$. Let $P_i$ be the edge of the geodesic loop which originates from $v_i$. Let $\omega_i$ be the unit vector at $v_i$ which is in direction of $P_{i-1}$ and $\omega'_i$ be the unit vector at $v_i$ which is in direction of $P_i$. Note that indices are given modulo $n$. See Figure \ref{polygonal-loop}.

\begin{lemma}
	\label{polygonal-holonomy}
	Let
	$$\Theta= \sum_{i=1}^{n}\theta(\Omega_i,\Omega'_{i}).$$
	We have
	$$hol(P)=e^{\imath \Theta}.$$
	\begin{proof}
		Let $\beta$ be a non-zero vector based at $V_i$ and
		$\beta'$ be the parallel transport of $\Omega$ through $P_i$.
		Since  parallel transport of  $\Omega'_i$ through $P_i$ is $\Omega_{i+1}$, we have
		
		$$\theta(\beta', \Omega'_{i+1})=\theta(\beta,\Omega'_{i})+\theta(\Omega_{i+1},\Omega'_{i+1}).$$
	Let $\alpha$ be a vector based at $V_1$. We denote its parallel transport at $V_k$ (through $P$) by $\alpha^{(k-1)}$. Above equality implies the followings:
	
	$$\theta(\alpha^{(1)},\Omega'_2)=\theta(\alpha,\Omega'_1)+\theta(\Omega_2,\Omega'_2)$$
	
		$$\theta(\alpha^{(2)},\Omega'_3)=\theta(\alpha^{(1)},\Omega'_2)+\theta(\Omega_3,\Omega'_3)$$
		$$\dots$$
	$$\theta(\alpha^{(n-1)},\Omega'_{n})=\theta(\alpha^{(n-2)},\Omega'_{n-1})+\theta(\Omega_n,\Omega'_n)$$
		$$\theta(\alpha^{(n)},\Omega'_{1})=\theta(\alpha^{(n-1)},\Omega'_{n})+\theta(\Omega_1,\Omega'_1)$$
	Note that we obtained the last equality by counting indices modulo $n$, and $\alpha^{(n)}$ is parallel transport of $\alpha$ through $P$. Adding above equalities we get
	
	$$\theta(\alpha^{(n)},\Omega'_1)=\theta(\alpha,\Omega'_1)+\Theta,$$
	which means that $hol(P)=e^{\imath \Theta}$.
	\end{proof}
\end{lemma}

\begin{proposition}
	\label{cover-holonomy}
	Let $\psi: S^* \rightarrow S$ be a flat covering,  $x \in S^*-\psi^{-1}(\frak{s} \cup \psi(\frak{s^*}))$ and $y=\psi(x)$, where $\frak{s}$ and $\frak{s^*}$ are  the sets of singular points of $S$ and $S^*$, respectively.
	
	\item
	The following diagram is commutative.
	
	$$\begin{array}[c]{ccc}
	\Pi_{1}(S^*-\psi^{-1}(\frak{s}\cup \psi(\frak{s^*})),x)	&\stackrel{hol}{\longrightarrow}& \mathbb{T}\\
	\Big\downarrow\scriptstyle{\psi_*}&&\Big\Updownarrow\scriptstyle{}\\
	\Pi_{1}(S-(\frak{s}\cup \psi(\frak{s^*})), y)&\stackrel{hol}{\longrightarrow}&\mathbb{T}
	\end{array}$$

	\begin{proof}
		
		Take a polygonal loop $P$ in $S^*-\psi^{-1}(\frak{s} \cup \psi(\frak{s^*}))$ which is based at $x$. By Lemma \ref{polygonal-holonomy},
		$hol(P)$ depends only sum of angles between edges of $P$. Observe that $\psi(P)$ is a polygonal loop in $S-(\frak{s}\cup \psi(\frak{s^*}))$ which is based at $\psi(x)$. Since $\psi$ is a local isometry, it is angle preserving. This
		implies that $hol(P)=hol(\psi_*(P))$.
		
	\end{proof}
\end{proposition}

Now we relate really flat surfaces to very flat surfaces.  Let
$S$ be a really flat surface. For each  $x \in S$, let
$\theta(x)=\frac{2\pi k(x)}{l(x)}$, where $k(x)$ and $l(x)$ are coprime positive integers. Let $\frak{s}$ be the set of singular points of $S$ and $y \in S- \frak{s}$. Let
$N$ be the order of the holonomy group of $S$ at $y$. Consider
$S$ as a differentiable orbifold with a divisor $D=\sum_{x \in \frak{s}}Nx$.

\begin{theorem}
	\label{criteria-really-flatness}
	\begin{enumerate}
		\item
		There exists a flat Galois covering $\psi: S^* \rightarrow S$ so that
		$S^*$ is very flat.
		
		\item
		Assume that $ \mathbf{\alpha}:  S''\rightarrow S'$ is a flat covering.  $S''$ is really flat if and only if $S'$ is really flat.
	\end{enumerate}
	\begin{proof}
		\begin{enumerate}
			\item
			For each $x\in S$, $N$ is divisible by $l(x)$. Therefore, there exists a homomorphism $\tilde{hol}$ so that the triangle on the right of  the following diagram is commutative:
			\[
			\begin{tikzcd}
			\Pi_{1}(S^*-\psi^{-1}(\frak{s}),y^*)\arrow{r}{\psi_*}\arrow{d}
			&	\Pi_{1}(S-\frak{s},y)\arrow{rd}{\text{hol}}\arrow{d}
			\\
			\Pi_{1}^{orb}(S^*,D^*,y^*)\arrow{r}{\psi_*}&\Pi_{1}^{orb}(S,D,y)\arrow{r}{\tilde{hol}}&
			\mathbb{T}
			\end{tikzcd}
			\]

			Consider $\frak{K}=\text{Ker}(\tilde{hol})$. Observe  that $\frak{K}$ is a finite index normal subgroup of $\Pi_{1}^{orb}(S,D,y)$. Consider  Galois orbifold covering  $\psi: (S^*,D^*) \rightarrow (S,D)$ which corresponds to $\frak{K}$. Observe that this covering is finite.  Let $y^*$ be a point  $S^*$ so that $\psi(y^*)=y$. Commutativity of the square in the above diagram  implies that the diagram is commutative. Consider the metric induced on $S^*$ from the one on $S$ by the map $\psi$. Observe that $\psi : S^* \rightarrow S$ is a flat map for the metrics on $S$ and $S^*$.
			Proposition \ref{cover-holonomy} implies that the following diagram is also commutative:
			\[
			\begin{tikzcd}
			\Pi_{1}(S^*-\psi^{-1}(\frak{s}),y^*)\arrow{r}{\psi_*} \arrow{rrd}{hol}
			&	\Pi_{1}(S-\frak{s},y)\arrow{rd}{\text{hol}}
			\\
			&&
			\mathbb{T}
			\end{tikzcd}
			\]
			
			Since $\tilde{hol}\circ \psi_* \ : \Pi_{1}^{orb}(S^*,D^*,y^*)\rightarrow \mathbb{T}$ is trivial homomorhism, commutativity of the first diagram implies that 
			$hol \circ \psi_*: \Pi_{1}(S^*-\psi^{-1}(\frak{s}),y^*)\rightarrow \mathbb{T}$ is also trivial. By the commutativity of above diagram,
			the map
			$$hol:	\Pi_{1}(S^*-\psi^{-1}(\frak{s}),y^*) \rightarrow \mathbb{T}$$
			is trivial, which means that $S^*$ is very flat.
			
			\item
			Assume that $S''$ is really flat. By the first item of the present theorem, we may assume that $S''$ is
			very flat. Let $\frak{b}$ be set of branched points of $\mathbf{\alpha}$. Let
			and $\frak{s'}, \ \frak{s''}$ be set of singular points of $S'$ and $S''$. Let $A=\frak{b}\cup \mathbf{\alpha}(\frak{s''})\cup \frak{s'}$. Let $y' \in S'- A$ and $y''\in \mathbf{\alpha}^{-1}(A)$. Consider the below diagram which is commutative by Proposition \ref{cover-holonomy}.
			
			\[
			\begin{tikzcd}
			\Pi_{1}(S''-\mathbf{\alpha}^{-1}(A),y'')\arrow{r}{\mathbf{\alpha}_*} \arrow{rrd}{hol''}
			&	\Pi_{1}(S'-A,y')\arrow{rd}{hol'}
			\\
			&&
			\mathbb{T}
			\end{tikzcd}
			\]
			Let $\frak{K}=Ker(hol')$. We have that $$\mathbf{\alpha}_*(\Pi_{1}(S''-\mathbf{\alpha}^{-1}(A),y'')) \subset \frak{K} \subset \Pi_{1}(S'-A,y').$$
			Observe  that
			$\lvert \Pi_{1}(S'-A,y') / \mathbf{\alpha}_*(\Pi_{1}(S''-\mathbf{\alpha}^{-1}(A),y''))   \rvert < \infty$, since $\mathbf{\alpha}$ is a finite covering. This implies that $\frak{K}$ is a  finite index subgroup of $\Pi_{1}(S'-A,y')$. Hence image of $hol'$ is finite. Thus, $S'$ is really flat.
			
			Other implication of the last part of the present proposition  immediately follows from Proposition \ref{cover-holonomy}.
		\end{enumerate}
	\end{proof}
\end{theorem}

\section{Ergodicity}

\ \ \ \ In this section, we study behavior of a typical geodesic and billiard trajectory on a really flat surface. Note that if the holonomy group of the surface is not trivial, we can not talk about directional flows. Therefore, we consider flow of  a generic geodesics on a really flat surfaces. First, we recall notion of  ergodicity,  and then state a main theorem which is proved by  Kerckhoff, Masur
and Smillie \cite{directional-ergodic}.

Let $X$ be a set, $\Sigma$ be a $\sigma-$algebra  and $\mu$ be a finite measure
on $X$.

\begin{udefinition}
	A measurable semiflow $T_t, \ t \in \mathbb{R}^+\cup \{0\}$ on  $(X,\Sigma, \mu)$ is called ergodic if for each measurable
	set $A$ of $X$ satisfying $\mu(T_t^{-1}(A)\ \Delta \ A)=0$ for all $t \in \mathbb{R}^+\cup \{0\}$,
	either $\mu(A)$ or $\mu(X-A)$ is equal to $0$.
\end{udefinition}

\begin{udefinition}
	A  measurable semiflow $T_t$ on a measurable space is called uniquely ergodic if there exists  a unique invariant measure  $\mu$ (up to scaling with a positive number) for which  $T_t$ is ergodic.
\end{udefinition}

\subsection{Ergodicity of geodesic flows on very flat surfaces}
\ \ \ \ Let $S'$ be a very flat surface. Observe that triviality of the holonomy group
implies that for each point on $x$, any vector on tangent space of $x$ can be carried to any non-singular point of $S'$ by a parallel transport and resulting vector is independent of the chosen path. That is, directional flow exists for any direction on the surface.

\begin{utheorem}[\bf{Kerchoff, Masur, Smillie}]
	\label{KMS}
	For almost all directions, directional flow on $S'$ is uniquely ergodic with respect to  the area measure of $S'$.
\end{utheorem}

Let's denote geodesic flow in a direction by $g_t$ and area measure of $S'$ by $\mu$. Note that we can state the result above as follows. For almost all
directions and for any continuous function $f$ on $S'$
\begin{align*}
\frac{1}{K}\int_{0}^{K} f(g_t(x)) dt \ \ \textup{converges uniformly to}\
\frac{1}{\mu(S')}\int_{S'} f d\mu.
\end{align*}

\subsection{Ergodicity of geodesic flows on  really flat surfaces without boundary}

\ \ \ \ \ Let $S$ be an orientable really flat surface without boundary. Let $\phi: \ S^* \rightarrow S$ be the flat Galois cover  which is constructed in Theorem \ref{criteria-really-flatness} and corresponds to the kernel of the homomorphism $hol \ : \Pi_1(S-\frak{s}) \rightarrow \mathbb{T}$. Let $C=\Pi_1(S -\frak{s})/Ker(hol)$ be the group of deck transformations of the cover $\phi: \ S^* \rightarrow S$. Note that $C$
acts transitively on the fiber of each point and $C$ is cyclic. Also observe that $C$ acts on $S^*$ by isometries. Let $N$ be the order of $C$, or equivalently, degree of the cover. Let $\mu$ and $\mu^*$ be the area measures on $S$ and $S^*$, respectively. We assume that $\mu(S)=1$, thus $\mu^*(S^*)=N$.

Now, we state a lemma whose proof is based on Theorem \ref{criteria-really-flatness}.

\begin{lemma}
	\label{main}
	For almost all $x \in S$ and for almost all $\omega \in T_x(S)$, flow of the geodesic in direction of $\omega$ is ergodic.
	\begin{proof}
		Consider the flat covering $\phi: S^* \rightarrow S$. Let $\frak{s}$ be the set of singular points of $S$. For each $x\in S-\frak{s}$, $\omega \in T_x(S)$, let $g_t(x)$ be the geodesic based at $x$ with initial direction $\omega$ and $g^*_t(x^*)$ be its lift  so that $g^*_0(x^*)=x^*$. Let $\omega^* \in T_{x^*}$ so that $d \phi_{x^*}(\omega^*)=\omega$. It is clear that initial direction of $g^*_t(x^*)$ is $\omega^*$.
		
		Observe that for almost all $x \in S$ and $\omega \in T_x(p)$, flow of the geodesic $g_t^*(x^*)$ is ergodic and $g_t^*(x^*)$ does not pass through $\phi^{-1}(\frak{s})$. See Section \ref{KMS}. Thus for each  such $x$ and $\omega$, the following the following limit exists
		
		\begin{align*}
		\lim_{K \rightarrow \infty}	\frac{1}{K}\int_{0}^{K} f(\phi(g_t^*(x^*))) dt \ \ =
		\frac{1}{N}	\int_{S^*} f\circ \phi \ d\mu^*  ,
		\end{align*}

		\noindent	for any continuous function $f$ on $S$. Since 
		
		$$f(\phi(g_t^*(x^*)))=f(g_t(x)) \textup{ and} \int_X f\circ \phi \ d\mu^*=N\int_S f \ d \mu,$$ the following limit also exists:
		
		\begin{align*}
		\lim_{K \rightarrow \infty}	\frac{1}{K}\int_{0}^{K} f(g_t(x)) dt \ \ =
		\int_{S} f \ d\mu.
		\end{align*}
		
	\end{proof}
\end{lemma}

We continue to use notation of the above proof. Assume that directional flow on
$S^*$ with respect to $\omega^*$ is ergodic. Let $(y,\omega_y) \in T(S)$ so that $\omega_y$ is obtained from parallel transport of $\omega$ through a path joining $x$ to $y$.

\begin{corollary}
	For almost all $y \in S$, flow of the geodesic $g_t(y)$ with initial direction $\omega_y$ is ergodic.
	\begin{proof}
		Consider directional  flow on $S^*$ obtained from $(x^*,\omega^*)$. For each $(y,\omega_y)$ having above properties, there is $(y^*,\omega_y^*) \in T(S^*)$ such that $\ y^* \in \phi^{-1}(y)$, $d_{y^*}(\omega_y^*)=\omega_y$ and $\omega_y^*$   is obtained by parallel transport of $\omega^*$.
		Since directional flow is ergodic, we see that for almost all $y$,
		geodesic flow with respect to $(y^*,\omega_y^*) \in T(S^*)$ is ergodic.
		As in the proof of the Lemma \ref{main}, this implies that for almost all $y$, geodesic flow with respect to $(y,\omega_y)$ is ergodic.
	\end{proof}
\end{corollary}

Now we generalize result of Lemma \ref{main} to non-orientable flat surfaces.

\begin{theorem}
	\label{geodesic-ergodic}
	Let $S'$ be a really flat surface without boundary. For almost all $x \in S'$ and for almost all $\omega \in T_x(S')$, flow of the geodesic in direction of $\omega$ is ergodic.
	\begin{proof}
		If $S'$ is orientable, then the statement is true. See Lemma \ref{main}. Assume that $S'$ is not orientable. Consider orientable double cover $\psi : \ S'' \rightarrow S'$. The map 
		$\psi$ induces a metric on $S''$ by the one on $S'$, and it is a flat map with respect to these metrics. Observe that flow of a geodesic on $S'$ is ergodic if and only if flow of one of its lifts is ergodic. Since $S''$ is really flat and orientable, Lemma \ref{main} implies the result.
	\end{proof}
\end{theorem}

\subsection{Ergodicity of billiard flows on really flat surfaces }

\ \ \ \ We use Theorem \ref{geodesic-ergodic}  to obtain a similar result about
ergodicity of billiard flows on really flat surfaces with boundary.

\begin{theorem}
	\label{ergodicity-billiard}
	Let $S$ be a really flat surface with boundary. For almost all $x \in S$ and for almost all $\omega \in T_x(S)$, flow of the billiard in direction of $\omega$ is ergodic.
	\begin{proof}
		Observe that flow of a billiard  on $S$ is ergodic if and only if flow of the corresponding
		geodesic in double of $S$, $S^D$, is ergodic. Since $S^D$ is really flat, result follows from Theorem  \ref{geodesic-ergodic}.
	\end{proof}
	
\end{theorem}

\label{references}


\end{document}